\documentclass[leqno,12pt]{article}

\usepackage[dvips]{graphics}
\usepackage{color}
\usepackage{ifthen}
\usepackage{soul}
\usepackage{latexsym}
\usepackage{amsthm,amssymb}
\usepackage{amsbsy,amsfonts,amsmath}

\usepackage{a4}

\numberwithin{equation}{section}

\newtheorem{theorem}{Theorem}[section]

\theoremstyle{definition}

\theoremstyle{remark}
\newtheorem{remark}[theorem]{Remark}

\newcommand{\rmSpan}{{\mathrm{Span}}}

\newcommand{\al}{\alpha}

\newcommand{\bN}{{\mathbb{N}}}
\newcommand{\bZ}{{\mathbb{Z}}}
\newcommand{\bZgeqo}{{\bZ_{\geq 0}}}

\newcommand{\bR}{{\mathbb{R}}}
\newcommand{\bQ}{{\mathbb{Q}}}
\newcommand{\bC}{{\mathbb{C}}}
\newcommand{\bK}{{\mathbb{K}}}
\newcommand{\bKt}{\bK^\times}

\newcommand{\rtminusone}{{\sqrt{-1}}}

\newcommand{\fkJ}{\mathbb{J}}

\newcommand{\fka}{{\mathfrak{a}}}
\newcommand{\fkgl}{{\mathrm{gl}}}
\newcommand{\fksl}{{\mathrm{sl}}}
\newcommand{\rmA}{{\mathrm{A}}}

\newcommand{\mate}{e}
\newcommand{\brvee}{I}
\newcommand{\brvpi}{{\breve \pi}}

\newcommand{\acutep}{{\acute p}}
\newcommand{\acuted}{{\acute d}}

\newcommand{\acuteh}{{\acute h}}
\newcommand{\acutee}{{\acute e}}
\newcommand{\acutef}{{\acute f}}

\newcommand{\str}{{\mathrm{str}}}

\newcommand{\checkc}{c}
\newcommand{\checkd}{{t{\frac {d} {dt}}}}

\newcommand{\defgamma}{\gamma}

\newcommand{\defg}{g}
\newcommand{\defphi}{{\varphi}}
\newcommand{\defpsi}{{\psi}}

\newcommand{\rmid}{{\mathrm{id}}}

\newcommand{\deftheta}{\theta}

\newcommand{\hatn}{m}

\newcommand{\deftau}{\tau}

\newcommand{\defnu}{n}

\newcommand{\defN}{N}
\newcommand{\defNp}{{\defN^\prime}}

\newcommand{\barp}{{\bar p}}
\newcommand{\bard}{{\bar d}}

\newcommand{\barh}{{\bar h}}
\newcommand{\bare}{{\bar e}}
\newcommand{\barf}{{\bar f}}

\newcommand{\defM}{M}
\newcommand{\hatdefN}{{\hat{\defN}}}
\newcommand{\bI}{{\mathbb{I}}}

\newcommand{\defhatdel}{{\hat \delta}}
\newcommand{\defprt}{\partial}
\newcommand{\bveps}{\bar{\varepsilon}}
\newcommand{\defV}{V}

\newcommand{\cU}{{\mathcal{U}}}
\newcommand{\fkntauo}{{\mathfrak{n}}^{(\deftau,0)}}
\newcommand{\fkntaup}{{\mathfrak{n}}^{(\deftau,+)}}
\newcommand{\fkntaum}{{\mathfrak{n}}^{(\deftau,-)}}

\newcommand{\defq}{q}

\newcommand{\hatdefV}{{\hat{\defV}}}
\newcommand{\hatdefVp}{\hatdefV^\prime}
\newcommand{\hatdefVprp}{\hatdefV^{\prime,+}}

\newcommand{\bhm}{\chi}
\newcommand{\hatbhm}{{\hat{\bhm}}}
\newcommand{\bhmop}{{\bhm^{\mathrm{op}}}}
\newcommand{\defU}{U}
\newcommand{\defUo}{\defU^0}
\newcommand{\defUp}{\defU^+}
\newcommand{\defUm}{\defU^-}
\newcommand{\defUbhm}{\defU_\bhm}
\newcommand{\defUobhm}{\defUo_\bhm}
\newcommand{\defUpbhm}{\defUp_\bhm}
\newcommand{\defUmbhm}{\defUm_\bhm}

\newcommand{\tildefU}{{\tilde{\defU}}}
\newcommand{\tildefUbhm}{\tildefU_\bhm}
\newcommand{\tildefUo}{{\tilde{\defU}}^0}
\newcommand{\tildefUp}{{\tilde{\defU}}^+}
\newcommand{\tildefUm}{{\tilde{\defU}}^-}
\newcommand{\tildefUobhm}{\tildefUo_\bhm}
\newcommand{\tildefUpbhm}{\tildefUp_\bhm}
\newcommand{\tildefUmbhm}{\tildefUm_\bhm}

\newcommand{\tilUpsbhmop}{{\tilde{\Upsilon}}_\bhmop}
\newcommand{\tildefUbhmop}{\tildefU_\bhmop}

\newcommand{\defK}{K}
\newcommand{\defL}{L}
\newcommand{\defE}{E}
\newcommand{\defF}{F}

\newcommand{\deftilK}{{\tilde{\defK}}}
\newcommand{\deftilL}{{\tilde{\defL}}}
\newcommand{\deftilE}{{\tilde{\defE}}}
\newcommand{\deftilF}{{\tilde{\defF}}}

\newcommand{\defDelta}{\Delta}
\newcommand{\defS}{S}
\newcommand{\defepsilon}{\epsilon}

\newcommand{\defvartheta}{\vartheta}
\newcommand{\tildefvartheta}{{\tilde{\defvartheta}}}
\newcommand{\tildefvarthetabhm}{\tildefvartheta_\bhm}

\newcommand{\tilIpbhm}{{\tilde{I}}^+_\bhm}
\newcommand{\tilImbhm}{{\tilde{I}}^-_\bhm}

\newcommand{\tilIpbhmop}{{\tilde{I}}^+_\bhmop}

\newcommand{\tilJbhm}{{\tilde{J}}_\bhm}

\newcommand{\tilpibhm}{{\tilde{\pi}}_\bhm}

\newcommand{\wbrl}{{[\![}}
\newcommand{\wbrr}{{]\!]}}

\newcommand{\tilIntpbhm}{{\tilde{I}}^{\natural,+}_\bhm}
\newcommand{\tilIntmbhm}{{\tilde{I}}^{\natural,-}_\bhm}
\newcommand{\tilIntpbhmop}{{\tilde{I}}^{\natural,+}_\bhmop}
\newcommand{\tilJntbhm}{{\tilde{J}}^\natural_\bhm}
\newcommand{\defUntbhm}{\defUbhm^\natural}
\newcommand{\defUntobhm}{\defUbhm^{0,\natural}}
\newcommand{\defUntpbhm}{\defUbhm^{+,\natural}}
\newcommand{\defUntmbhm}{\defUbhm^{-,\natural}}
\newcommand{\tilpintbhm}{{^\natural\tilpibhm}}

\newcommand{\defntK}{\defK^\natural}
\newcommand{\defntL}{\defL^\natural}
\newcommand{\defntE}{\defE^\natural}
\newcommand{\defntF}{\defF^\natural}

\newcommand{\pintbhm}{\pi^\natural}

\newcommand{\defPsi}{\Psi}

\newcommand{\defntA}{{\mathcal{A}}^\natural}
\newcommand{\defntB}{{\mathcal{B}}^\natural}

\newcommand{\defntZ}{{\mathcal{Z}}^\natural}

\newcommand{\defs}{s}
\newcommand{\defa}{a}

\newcommand{\defsigma}{\sigma}
\newcommand{\defwp}{\wp}
\newcommand{\hatdefwp}{{\hat{\defwp}}}

\newcommand{\dl}{\defhatdel}
\newcommand{\be}{\bveps}

\newcommand{\DiagramAA}{
\setlength{\unitlength}{1mm}
\begin{picture}(140,60)(7,5)

\put(10,5){$(\deftau=1,\,\,\fksl^{(1)}_\barp=\fksl^{(1)}(\defnu|\defnu),\,\,\defnu={\frac \defN 2},\,\,\defN\geq 4)$}

\put(70, 55){\circle{5}}
\put(68.1, 53.1){\line(1,1){3.7}} \put(69, 60){$0$}
\put(62,45){$\dl-\be_1+\be_\defN$}
\put(68.1, 56.9){\line(1,-1){3.7}}
\put(72.5, 55){\line(1,0){48}}
\put(67.5, 55){\line(-1,0){48}}

\put(19.5,40){\oval(20,30)[l]}

\put(22, 25){\circle{5}}
\put(21, 30){$1$}
\put(16,15){$\be_1-\be_2$}
\put(24.5, 25){\line(1,0){7}}

\put(35, 24.7){\dots}

\put(43.5, 25){\line(1,0){7}}
\put(53, 25){\circle{5}}
\put(49, 30){$\defnu-1$}
\put(40,15){$\be_{\defnu-2}-\be_{\defnu-1}$}
\put(55.5, 25){\line(1,0){12}}

\put(70, 25){\circle{5}}
\put(68.1, 23.1){\line(1,1){3.7}} \put(69, 30){$\defnu$}
\put(62,15){$\be_{\defnu-1}-\be_\defnu$}
\put(68.1, 26.9){\line(1,-1){3.7}}
\put(72.5, 25){\line(1,0){7}}

\put(77.5, 25){\line(1,0){7}}
\put(87, 25){\circle{5}}
\put(83, 30){$\defnu+1$}
\put(80,15){$\be_{\defnu}-\be_{\defnu+1}$}
\put(89.5, 25){\line(1,0){7}}

\put(100, 24.7){\dots}

\put(108.5, 25){\line(1,0){7}}
\put(118, 25){\circle{5}}
\put(112, 30){$\defN-1$}
\put(108,15){$\be_{\defN-1}-\be_{\defN}$}

\put(120.5,40){\oval(20,30)[r]}

\end{picture}}

\newcommand{\DiagramCC}{
\setlength{\unitlength}{1mm}
\begin{picture}(140,30)(-5,-5)

\put(0,5){$(\deftau=2,\,\,\fksl^{(2)}_\barp=\fksl^{(2)}(\defN|\defN),\,\,\defnu={\frac \defN 2},\,\,\defN\geq 4)$}

\put(05, 25){\circle{5}}
\put(04, 30){0} \put(00,15){$\dl-2\be_1$} 
\put(19.5, 25){\line(-1,-1){3.7}}
\put(19.5, 25){\line(-1,1){3.7}}
\put(18, 23.5){\line(-1,0){11}}
\put(18, 26.5){\line(-1,0){11}}

\put(22, 25){\circle{5}}
\put(21, 30){$1$}
\put(16,15){$\be_1-\be_2$}
\put(24.5, 25){\line(1,0){2}}

\put(30, 24.7){\dots}

\put(38.5, 25){\line(1,0){2}}
\put(43, 25){\circle{5}}
\put(39, 30){$\defnu-1$}
\put(30,15){$\be_{\defnu-2}-\be_{\defnu-1}$}
\put(45.5, 25){\line(1,0){12}}

\put(60, 25){\circle{5}}
\put(58.1, 23.1){\line(1,1){3.7}} \put(59, 30){$\defnu$}
\put(52,15){$\be_{\defnu-1}-\be_\defnu$}
\put(58.1, 26.9){\line(1,-1){3.7}}
\put(62.5, 25){\line(1,0){7}}

\put(67.5, 25){\line(1,0){7}}
\put(77, 25){\circle{5}}
\put(73, 30){$\defnu+1$}
\put(70,15){$\be_{\defnu}-\be_{\defnu+1}$}
\put(79.5, 25){\line(1,0){2}}

\put(85, 24.7){\dots}

\put(93.5, 25){\line(1,0){2}}
\put(98, 25){\circle{5}}
\put(92, 30){$\defN-1$}
\put(88,15){$\be_{\defN-1}-\be_{\defN}$}
\put(100.5, 25){\line(1,1){3.7}}
\put(100.5, 25){\line(1,-1){3.7}}
\put(102, 23.5){\line(1,0){11}}
\put(102, 26.5){\line(1,0){11}}

\put(115, 25){\circle{5}}
 \put(114, 30){$\defN$} \put(113,15){$2\be_{\defN}$}

\end{picture}}

\newcommand{\DiagramBB}{
\setlength{\unitlength}{1mm}
\begin{picture}(140,30)(-5,-5)

\put(0,5){$(\deftau=4,\,\,\fksl^{(4)}_\barp=\fksl^{(4)}(\defN+1|\defN+1),\,\,\defnu={\frac \defN 2},\,\,\defN\geq 2)$}
\put(05, 25){\circle{5}}
\put(04, 30){0} \put(00,15){$\dl-\be_1$} 
\put(07.5, 25){\line(1,1){3.7}}
\put(07.5, 25){\line(1,-1){3.7}}
\put(09, 23.5){\line(1,0){11}}
\put(09, 26.5){\line(1,0){11}}

\put(22, 25){\circle{5}}
\put(21, 30){$1$}
\put(16,15){$\be_1-\be_2$}
\put(24.5, 25){\line(1,0){2}}

\put(30, 24.7){\dots}

\put(38.5, 25){\line(1,0){2}}
\put(43, 25){\circle{5}}
\put(39, 30){$\defnu-1$}
\put(30,15){$\be_{\defnu-2}-\be_{\defnu-1}$}
\put(45.5, 25){\line(1,0){12}}

\put(60, 25){\circle{5}}
\put(58.1, 23.1){\line(1,1){3.7}} \put(59, 30){$\defnu$}
\put(52,15){$\be_{\defnu-1}-\be_\defnu$}
\put(58.1, 26.9){\line(1,-1){3.7}}
\put(62.5, 25){\line(1,0){7}}

\put(67.5, 25){\line(1,0){7}}
\put(77, 25){\circle{5}}
\put(73, 30){$\defnu+1$}
\put(70,15){$\be_{\defnu}-\be_{\defnu+1}$}
\put(79.5, 25){\line(1,0){2}}

\put(85, 24.7){\dots}

\put(93.5, 25){\line(1,0){2}}
\put(98, 25){\circle{5}}
\put(92, 30){$\defN-1$}
\put(88,15){$\be_{\defN-1}-\be_{\defN}$}

\put(115, 25){\circle{5}}
 \put(114, 30){$\defN$} \put(114,15){$\be_{\defN}$} 
\put(112.5, 25){\line(-1,-1){3.7}}
\put(112.5, 25){\line(-1,1){3.7}}
\put(111, 23.5){\line(-1,0){11}}
\put(111, 26.5){\line(-1,0){11}}

\end{picture}}

\begin{document}

\begin{center}
{\Large{
Lowest positive almost central elements of \\ $U_q(\fksl^{(1)}(n|n))$ $(n\geq 2)$, 
$U_q(\fksl^{(2)}(2n|2n))$ $(n\geq 2)$ 
and $U_q(\fksl^{(4)}(2n+1|2n+1))$ $(n\geq 1)$ and \\ their multi-parameter 
quantum affine superalgebras}}
\\ \quad \\ 
Hiroyuki Yamane
\end{center}

\begin{abstract}
We study $U_qsl^{(\tau)}(M|M)$, where $\tau$ is $1$, $2$ or $4$.
Let $\grave{\pi}:sl^{(\tau)}(M|M)\to A^{(\tau)}(M-1,M-1)$ be the canonical projection.
Let $\breve{e}$ be the $2M\times 2M$ unit matrix.
Then $\ker \grave{\pi}$
$=\mathbb{C}[t^s,t^{-s}]\otimes \breve{e}(=\oplus_{k=-\infty}^\infty\mathbb{C} t^{ks}\otimes \breve{e})$,
where $s=1$ (resp. $=2$) if $\tau$ is $1$ or $2$ (resp. is $4$).
In this paper, we give an explicit formula of the element $\defntZ$
(see \eqref{eqn:defZnatural} and \eqref{eqn:profdefZnatural}) of
$U_qsl^{(\tau)}(M|M)$ corresponding to
$t^s \otimes \breve{e}$, and do the same for its multi-parameter versions.
\end{abstract}

\section{Introduction }
Let $m$, $n\in\bN$.
Let $\Im_{m|n}:=\{k\in\bZ|1\leq k\leq m+n\}$,
$\Im_{m|n}^\prime:=\{k\in\bZ|1\leq k\leq m\}$ and
$\Im_{m|n}^{\prime\prime}:=\{k\in\bZ|m+1\leq k\leq m+n\}$.
Denote by ${\mathrm{M}}_{m+n}(\bC)$ 
the associative $\bC$-algebra formed by
all the $(m+n)\times (m+n)$-matrices over $\bC$.
For $i$, $j\in\Im_{m|n}$, let $\mate_{ij}\in{\mathrm{M}}_{m+n}(\bC)$ denote
the $(i,j)$-matrix unit, so $\{\mate_{ij}|i,j\in\Im_{m|n}\}$ is the $\bC$-basis of ${\mathrm{M}}_{m+n}(\bC)$
with $\mate_{ij}\mate_{kl}=\delta_{jk}\mate_{il}$,
where $\delta_{jk}$ denotes the Kronecker's delta.
The Lie superalgebra $\fkgl(m|n)$ is such that 
it is an $(m+n)^2$-dimensional $\bC$-linear space identified with
${\mathrm{M}}_{m+n}(\bC)$, its even 
(resp. odd) part $\fkgl(m|n)(0)$ (resp. $\fkgl(m|n)(1)$) is its 
$m^2+n^2$-dimensional (resp. $2mn$-dimensional) subspace with 
the $\bC$-basis $\{\mate_{ij}|(i,j)\in(\Im_{m|n}^\prime)^2\cup(\Im_{m|n}^{\prime\prime})^2\}$
(resp. $\{\mate_{ij}|(i,j)\in(\Im_{m|n}^\prime\times\Im_{m|n}^{\prime\prime})\cup
(\Im_{m|n}^{\prime\prime}\times\Im_{m|n}^\prime)\}$),
and its super-bracket is defined by $[X,Y]:=XY-(-1)^{st}YX$ 
($X\in\fkgl(m|n)(s)$, $Y\in\fkgl(m|n)(t)$).
The Lie superalgebra $\fksl(m|n)$ is the sub-Lie superalgebra of $\fkgl(m|n)$
with the $\bC$-basis $\{\mate_{ij}|i,j\in\Im_{m|n}, i\ne j\}
\cup\{\mate_{ii}-(-1)^{\delta_{im}}\mate_{i+1i+1}|i\in\Im_{m|n}\setminus\{m+n\}\}$.

Let $\brvee_{m|n}:=\sum_{i=1}^m\mate_{ii}+{\frac m n}\sum_{j=m+1}^{m+n}\mate_{jj}\in\fksl(m|n)$.

Here assume $m=n$. Since $\brvee_{n|n}$ is the unit matrix,
$\bC\brvee_{n|n}$ is an one dimensional ideal of $\fksl(n|n)$.
Let $\fksl^{(1)}(n|n):=\fksl(n|n)\otimes\bC[t,t^{-1}]\oplus\bC\checkc\oplus\bC\checkd$
be the affine-type Lie superalgebra of $\fksl(n|n)$.
Let $Z^{(1)}_k:=\brvee_{n|n}\otimes t^k\in\fksl^{(1)}(n|n)$ ($k\in\bZ$).
Then $\oplus_{k=-\infty}^\infty\bC Z^{(1)}_k$ is an ideal of $\fksl^{(1)}(n|n)$.
We have: 
\begin{equation}\label{eqn:prZk}
\begin{array}{l}
\mbox{(i) $[X\otimes t^r,Z^{(1)}_k]=0$ ($X\in\fksl(n|n)$, $r\in\bZ$),
$[\checkc,Z^{(1)}_k]=0$, $[\checkd,Z^{(1)}_k]=kZ^{(1)}_k$} \\
\mbox{\quad\, for $k\in\bZ$.} \\
\mbox{(ii) $[Z^{(1)}_k,Z^{(1)}_s]=0$ for $k$, $s\in\bZ$.}
\end{array}
\end{equation}

{\it{We inspire that the equations of \eqref{eqn:prZk} would become applied to
discover important unknown theories of physics or mathematics.}}

Let $Z^{(2)}_k:=Z^{(1)}_{2k-1}$ and $Z^{(4)}_k:=Z^{(1)}_{4k-2}$
for $k\in\bZ$.
If $n\geq 2$ (resp. $n\geq 3$),
we have the twisted affine-type Lie superalgebra $\fksl^{(2)}(n|n)$
(resp, $\fksl^{(4)}(n|n)$), which can be realized as a sub-Lie superalgebra of 
$\fksl^{(1)}(n|n)$.
Then $\oplus_{k=-\infty}^\infty\bC Z^{(\deftau)}_k$ is an ideal of $\fksl^{(\deftau)}(n|n)$ 
for $\deftau\in\{1,2,4\}$.

In this paper, we explicitly give an element \eqref{eqn:defZnatural} of 
$U_q(\fksl^{(\deftau)}(n|n))$ ($\deftau\in\{1,2,4\}$) or their multi-parameter version
corresponding to $Z^{(\deftau)}_1$ (i.e., only for $k=1$),
and we show that it satisfies equations \eqref{eqn:profdefZnatural}
similar to those of (i) of \eqref{eqn:prZk}.

As for the elements of $U_q(\fksl^{(\deftau)}(n|n))$ ($\deftau\in\{1,2,4\}$) 
corresponding to $Z^{(1)}_k$ for all $k\in\bZ$, 
they have already been given by \cite{Y99} for $\tau=1$,
and their existence for $\deftau\in\{1,2,4\}$ has been shown by \cite{EG09}.
\newline
\section{Lie superalgebras}
\subsection{The affine Lie superalgebras $A^{(1)}(m,n)$, $A^{(2)}(m,n)$, \newline $A^{(4)}(m,n)$}
\label{subsection:defAtaumn}
For $x$, $y\in\bR$, let
$\fkJ_{x,y}:=\{n\in\bZ|x\leq n\leq y\}$
and $\fkJ_{x,\infty}:=\{n\in\bZ|x\leq n\}$.

Let $\bK$ be a field of characteristic $0$.
Assume that there exists an element 
$\rtminusone$ of $\bK$ such that 
$\rtminusone^2=-1$.
Let 
$\fka$ be a $\bK$-linear space.
Suppose that $\fka$ is a direct sum of its subspaces $\fka(0)$ and $\fka(1)$,
so $\fka=\fka(0)\oplus\fka(1)$.
Let $\fka(2n):=\fka(0)$ and $\fka(2n+1):=\fka(1)$ for $n\in\bZ$.
Let $[\,,\,]:\fka\times\fka\to\fka$ be a $\bK$-bilinear map.
We call $\fka=(\fka,[\,,\,])$
{\it{a Lie superalgebra}} if 
\begin{equation*}
\begin{array}{l}
\mbox{$[X,Y]\in\fka(a+b)\quad (X\in\fka(a),Y\in\fka(b))$,} \\
\mbox{$[X,Y]=-(-1)^{ab}[Y,X]\quad (X\in\fka(a),Y\in\fka(b))$,} \\
\mbox{$[X,[Y,Z]]=[[X,Y],Z]+(-1)^{ab}[Y,[X,Z]]\quad (X\in\fka(a),Y\in\fka(b),Z\in\fka)$.}
\end{array}
\end{equation*}
For a Lie superalgebra $\fka$, 
let $\cU(\fka)$ denote the universal enveloping superalgebra of $\fka$. 

Let $n\in\bN$.
Let $\acutep:\fkJ_{1,n}\to\fkJ_{0,1}$ be a map.
Let $n_a:=|\acutep^{-1}(a)|$, so $n=n_0+n_1$.
Assume $n_0\geq 1$ and $n_1\geq 1$.
Let $\fkgl_\acutep$ be an $n^2$-dimensional $\bK$-linear space
with a basis $\{\mate_{ij}|i,j\in\fkJ_{1,n}\}$. 
We regard $\fkgl_\acutep$ as the Lie superalgera
by
\begin{equation*}
\begin{array}{l}
\fkgl_\acutep(0)
=\oplus_{\acutep(i)+\acutep(j)\in 2\bZ}\bK\mate_{ij},\\
\fkgl_\acutep(1)
=\oplus_{\acutep(i)+\acutep(j)\in 2\bZ+1}\bK\mate_{ij}, \\
\mbox{$[\mate_{ij},\mate_{kl}]=\delta_{jk}\mate_{il}
-(-1)^{(\acutep(i)+\acutep(j))(\acutep(k)+\acutep(l))}\delta_{il}\mate_{kj}$.}
\end{array}
\end{equation*}
Note that $[X,Y]=XY-(-1)^{ab}YX$ holds for $X\in\fkgl_\acutep(a)$
and  $Y\in\fkgl_\acutep(b)$, where $XY$ and $YX$ means the matrix-multiplication.
Define the $\bK$-linear map $\str_\acutep:\fkgl_\acutep\to\bK$ by
$\str_\acutep(\mate_{ij}):=\delta_{ij}(-1)^{\acutep(i)}$.
Then $\str_\acutep$ is a Lie superalgebra homomorphism.
Let $\fksl_\acutep:=\ker \str_\acutep$.
Then $\fksl_\acutep$ is a sub-Lie superalgebra of $\fkgl_\acutep$.
Note that $\fkgl_\acutep$ and $\fksl_\acutep$ are isomorphic to
$\fkgl_{\acutep^\prime}$ and $\fksl_{\acutep^\prime}$
for a map $\acutep^\prime:\fkJ_{1,n}\to\fkJ_{0,1}$
with $n_0=|(\acutep^\prime)^{-1}(0)|$. Let 
\begin{equation*}
\acuted_i:=(-1)^{\acutep(i)}\quad (i\in\fkJ_{1,n}), 
\end{equation*} and 
\begin{equation}\label{eqn:genofslnzeronone}
\begin{array}{l}
\acuteh_i:=\acuted_i\mate_{i,i}-\acuted_{i+1}\mate_{i+1,i+1},\,\,\acutee_i:=\acuted_i\mate_{i,i+1},\,\,
\acutef_i:=\acuted_i\mate_{i+1,i} \\
(i\in\fkJ_{1,n-1}).
\end{array}
\end{equation} We have
\begin{equation}\label{eqn:genofslhef}
\begin{array}{l}
\mbox{$[\acuteh_i,\acuteh_j]=0,\,\,[\acuteh_i,\acutee_j]=\str_\acutep(\acuteh_i\acuteh_j)\acutee_j,\,\,
[\acuteh_i,\acutef_j]=-\str_\acutep(\acuteh_i\acuteh_j)\acutef_j$}, \\
\mbox{$[\acutee_i,\acutef_j]=\delta_{ij}\acuteh_i$}, \\
\mbox{$[\acutee_i,\acutee_j]=[\acutef_i,\acutef_j]=0$ ($\str_\acutep(\acuteh_i\acuteh_j)=0$)}, \\
\mbox{$[\acutee_i,[\acutee_i,\acutee_j]]=[\acutef_i,[\acutef_i,\acutef_j]]=0$ ($\str_\acutep(\acuteh_i^2)\ne 0$,
$i\ne j$, $\str_\acutep(\acuteh_i\acuteh_j)\ne 0$)}, \\
\mbox{$[\acutee_j,[[\acutee_i,\acutee_j],\acutee_k]]=[\acutef_j,[[\acutef_i,\acutef_j],\acutef_k]]=0$ ($\str_\acutep(\acuteh_j^2)=0$,
$i\ne j\ne k\ne i$,} \\
\mbox{$\str_\acutep(\acuteh_i\acuteh_j)\ne 0$, $\str_\acutep(\acuteh_i\acuteh_k)\ne 0$, $\str_\acutep(\acuteh_j\acuteh_k)=0$)}.
\end{array}
\end{equation}
The elements of \eqref{eqn:genofslnzeronone} generate $\fksl_\acutep$ as a Lie superalgebra.
It is well-known (see \cite{Y94}) that
\begin{equation*}
\begin{array}{l}
\mbox{as a Lie superalgebra, $\fksl_\acutep$ is presented by the generators \eqref{eqn:genofslnzeronone}}\\
\mbox{and the relations \eqref{eqn:genofslhef}.} \\
\end{array}
\end{equation*}
The Lie superalgebras $\fkgl_\acutep$ and $\fksl_\acutep$ are denoted by
$\fkgl(n_0|n_1)$ and $\fksl(n_0|n_1)$ respectively.
If $n_0\ne n_1$,  $\fksl(n_0|n_1)$ is also denoted by $\rmA(n_0-1,n_1-1)$.
Let $\brvee:=\sum_{i=1}^n\mate_{ii}$, so $\brvee$ is the unit matrix. 
Note that $\brvee\in\fksl_\acutep$ if and only if $n_0=n_1$,
and, if this is the case, we have the Lie supleralgebra $\fksl_\acutep/\bK\brvee$,
which is denoted by $\rmA(n_0-1,n_1-1)$.
To emphasis $\acutep$, we also denote $\rmA(n_0-1,n_1-1)$ by $\rmA_\acutep$, that is, we mean
$\rmA_\acutep$ to be $\fksl_\acutep$ (resp. $\fksl_\acutep/\bK\brvee$)
if  $n_0\ne n_1$ (resp. $n_0=n_1$).
Define the Lie superalgebra isomorphism (resp. eimomorphism) 
$\brvpi:\fksl_\acutep\to\rmA_\acutep$ by
$\brvpi(X):=X$ (resp. $X+\bK\brvee$) for $X\in\fksl_\acutep$
if  $n_0\ne n_1$ (resp. $n_0=n_1$).

Note that $\str_\acutep(YX)=(-1)^{ab}\str_\acutep(XY)$ holds for $X\in\fkgl_\acutep(a)$
and  $Y\in\fkgl_\acutep(b)$. Then we see that
\begin{equation*}
\str_\acutep([X,Y]Z)=\str_\acutep(X[Y,Z])\quad (X,Y,Z\in\fkgl_\acutep).
\end{equation*}
We define the Lie superalgera  $\fkgl^{(1)}_\acutep$ by
\begin{equation*}
\begin{array}{l}
\fkgl^{(1)}_\acutep=\fkgl_\acutep\otimes\bK[t,t^{-1}]\oplus\bK\checkc\oplus\bK\checkd, \\
\mbox{$[X\otimes t^x,Y\otimes t^y]:=[X,Y]\otimes t^{x+y}+m\delta_{x+y,0}\str_\acutep(XY)\checkc$}, \\
\mbox{$[\checkc,X\otimes t^x]=[\checkc,\checkd]=0$}, \quad
\mbox{$[\checkd,X\otimes t^x]=xX\otimes t^x$}, 
\end{array}
\end{equation*} where $X,Y\in\fkgl_\acutep$ and $x,y\in\bZ$.
Define the Lie superalgebras
\begin{equation*}
\begin{array}{l}
\fksl^{(1)}_\acutep:=\fksl_\acutep\otimes\bK[t,t^{-1}]\oplus\bK\checkc\oplus\bK\checkd, \\
\rmA^{(1)}_\acutep:=\rmA_\acutep\otimes\bK[t,t^{-1}]\oplus\bK\checkc\oplus\bK\checkd 
\end{array}
\end{equation*} such that
$\fksl^{(1)}_\acutep$ is a sub-Lie superalgebra
of $\fkgl^{(1)}_\acutep$ and
we have the Lie superalgebra homomorphism
$\brvpi^{(1)}:\fksl^{(1)}_\acutep\to\rmA^{(1)}_\acutep$ define by
$\brvpi^{(1)}(X\otimes t^x):=\brvpi(X)\otimes t^x$,
$\brvpi^{(1)}(\checkc)=\checkc$ and $\brvpi^{(1)}(\checkd)=\checkd$.
We also denote $\fkgl^{(1)}_\acutep$, $\fksl^{(1)}_\acutep$, $\rmA^{(1)}_\acutep$ by
$\fkgl^{(1)}(n_1|n_2)$, $\fksl^{(1)}(n_1|n_2)$, $\rmA^{(1)}(n_1-1,n_2-1)$ respectively.

For $x\in\bN$, define the map $\defgamma_x:\fkJ_{1,x}\to\fkJ_{1,x}$ by $\defgamma_x(i):=x+1-i$.

Let $n\in\bN$.
Let $\acutep:\fkJ_{1,n}\to\fkJ_{0,1}$ be a map.
Assume that 
$\acutep(\defgamma_n(i))=\acutep(i)$ for all $i\in\fkJ_{1,n}$.
Define a map $\defg^\acutep:\fkJ_{1,n}\to\{1,-1\}$ by $\defg^\acutep(i)\defg^\acutep(\defgamma_n(i))=(-1)^{\acutep(i)}$
($i\in\fkJ_{1,n}$) and
$\defg^\acutep(j)=1$ ($j\in\fkJ_{1,{\frac {n+1} 2}}$).
Then we can define the Lie superalgebra automorphism 
$\defphi_\acutep:\fkgl_\acutep\to\fkgl_\acutep$ by
\begin{equation*}
\defphi_\acutep(\mate_{ij}):=-(-1)^{\acutep(j)(\acutep(i)+\acutep(j))}\defg^\acutep(i)\defg^\acutep(j)
\mate_{\defgamma_n(j)\defgamma_n(i)}
\quad(i,j\in\fkJ_{1,n}).
\end{equation*} We have
\begin{equation*}
\defphi_\acutep^2=\rmid_{\fkgl_\acutep},\,\,\defphi_\acutep(\fksl_\acutep)=\fksl_\acutep,\,\,
\defphi_\acutep(\brvee)=-\brvee.
\end{equation*}
Then we have a sub-Lie superalgebra $\fkgl^{(2)}_\acutep$ of $\fkgl^{(1)}_\acutep$ defined by
\begin{equation*}
\fkgl^{(2)}_\acutep:=(\bigoplus_{x=-\infty}^\infty\ker((-1)^x\rmid_{\fkgl_\acutep}-\defphi_\acutep)\otimes t^x)
\oplus\bK\checkc\oplus\bK\checkd.
\end{equation*} Let
$\fksl^{(2)}_\acutep:=\fksl^{(1)}_\acutep\cap\fkgl^{(2)}_\acutep$
and $\rmA^{(2)}_\acutep:=\brvpi^{(1)}(\fksl^{(2)}_\acutep)$.
We also denote $\fkgl^{(2)}_\acutep$, $\fksl^{(2)}_\acutep$, $\rmA^{(2)}_\acutep$ by
$\fkgl^{(2)}(n_1|n_2)$, $\fksl^{(2)}(n_1|n_2)$, $\rmA^{(2)}(n_1-1,n_2-1)$ respectively.

Let $n\in\bN$. Assume $n\in2\bN$ and $n\geq 4$. 
Let $\acutep:\fkJ_{1,n}\to\fkJ_{0,1}$ be a map. Let $\hatn:={\frac {n-2} 2}$.
Define  $\deftheta:\fkJ_{0,n-1}\to\fkJ_{1,n}$ by $\deftheta(i):=i+1$.
Let $\acutep:\fkJ_{1,n}\to\fkJ_{0,1}$ be a map.
Assume that $\acutep(0)=1$,  $\acutep(\deftheta(\hatn+1))=0$ and 
$\acutep(\deftheta(i))=\deftheta(\defgamma_{n-1}(i))$
for $i\in\fkJ_{1,\hatn+1}$.
Define the map $\defg^{\prime,\acutep}:\fkJ_{1,n-1}$ 
 by $\defg^{\prime,\acutep}(i)\defg^{\prime,\acutep}(\defgamma_{n-1}(i))=(-1)^{\acutep(\deftheta(i))}$
($i\in\fkJ_{1,n-1}$) and
$\defg(j)=1$ ($j\in\fkJ_{1,\hatn+1}$).
Then we can define the Lie superalgebra automorphism 
$\defpsi_\acutep:\fkgl_\acutep\to\fkgl_\acutep$ by
\begin{equation*}\begin{array}{l}
\defpsi_\acutep(\mate_{\deftheta(i)\deftheta(j)}) \\
\,\,=
\left\{\begin{array}{l}
-\mate_{\deftheta(0)\deftheta(0)} \quad\quad\quad\mbox{if $\deftheta(i)=\deftheta(j)=\deftheta(0)$,} \\
-\rtminusone\defg^{\prime,\acutep}(\defgamma_{n-1}(i))\mate_{\deftheta(0)\deftheta(\defgamma_{n-1}(i))} \quad\quad
\quad\mbox{if $\deftheta(i)\ne \deftheta(0)$ and $\deftheta(j)=\deftheta(0)$,} \\
-\rtminusone\defg^{\prime,\acutep}(j)\mate_{\deftheta(\defgamma_{n-1}(j))\deftheta(0)} \quad\quad
\quad\mbox{if $\deftheta(i)=\deftheta(0)$ and $\deftheta(j)\ne \deftheta(0)$,} \\
-(-1)^{\acutep(\deftheta(j))(\acutep(\deftheta(i))+\acutep(\deftheta(j))))}
\defg^{\prime,\acutep}(i)\defg^{\prime,\acutep}(j)\mate_{\deftheta(\defgamma_{n-1}(j))\deftheta(\defgamma_{n-1}(i))} \\
\quad\quad\quad\mbox{if $\deftheta(i)\ne\deftheta(0)$ and $\deftheta(j)\ne\deftheta(0)$.}
\end{array}\right.
\end{array}
\end{equation*}

We have
\begin{equation*}
\defpsi_\acutep^2\ne\rmid_{\fkgl_\acutep},\,\,\defpsi_\acutep^4=\rmid_{\fkgl_\acutep},\,\,\defpsi_\acutep(\fksl_\acutep)=
\fksl_\acutep,\,\,
\defpsi_\acutep(\brvee)=-\brvee.
\end{equation*}
Then we have a sub-Lie superalgebra $\fkgl^{(4)}_\acutep$ of $\fkgl^{(1)}_\acutep$ defined by
\begin{equation*}
\fkgl^{(4)}_\acutep:=(\bigoplus_{x=-\infty}^\infty\ker((\rtminusone)^x\rmid_{\fkgl_\acutep}-\defpsi_\acutep)\otimes t^x)
\oplus\bK\checkc\oplus\bK\checkd.
\end{equation*} Let
$\fksl^{(4)}_\acutep:=\fksl^{(1)}_\acutep\cap\fkgl^{(4)}_\acutep$
and $\rmA^{(4)}_\acutep:=\brvpi^{(1)}(\fksl^{(4)}_\acutep)$.
We also denote $\fkgl^{(4)}_\acutep$, $\fksl^{(4)}_\acutep$, $\rmA^{(4)}_\acutep$ by
$\fkgl^{(4)}(n_1|n_2)$, $\fksl^{(4)}(n_1|n_2)$, $\rmA^{(4)}(n_1-1,n_2-1)$ respectively.

\subsection{The Lie superagerbas $\fksl^{(\deftau)}(\defM|\defM)$ with $\deftau\in\{1,2,4\}$, \newline
where 
$\defM:={\frac \defN 2}$ (resp. $\defN$, resp. $\defN+1$) if $\deftau=1$
(resp. $2$, resp. $4$).}
\label{subsection:sltau}

We use the notation as follows.
As defined below, we will let the meaning of the symbol $\defnu$ be different than that used above. 
\begin{equation*}
\begin{array}{l}
\mbox{From now on until end of this paper,} \\ 
\mbox{let 
$\deftau\in\{1,2,4\}(=\fkJ_{1,4}\setminus\{2\})$, $\defnu\in\bN$, $\defN:=2\defnu$,} \\
\mbox{assume $\defnu\geq 2$ (resp. $\defnu\geq 1$)
if $\deftau\in\{1,2\}$ (resp. $\deftau=4$),} \\
\mbox{let $\hatdefN:=\defN-1$ (resp. $\defN$) if $\deftau=1$ (resp. $\deftau\in\{2,4\}$), let
$\bI:=\fkJ_{0,\hatdefN}$,} \\
\mbox{let $\defNp:=\defN$ (resp. $2\defN$, resp. $2\defN+2$)
if $\deftau=1$ (resp. $2$, resp. $4$),} \\
\mbox{and let $\defM:={\frac \defN 2}$ (resp. $\defN$, resp. $\defN+1$) if $\deftau=1$
(resp. $2$, resp. $4$).}
\end{array}
\end{equation*}
Define the map $\barp:\fkJ_{1,\defNp}\to\fkJ_{0,1}$ by $\barp(i):=0$ ($i\in\fkJ_{1,\defnu}$
(resp. $\fkJ_{1,\defnu}\cup\fkJ_{3\defnu+1,4\defnu}$, 
resp. $\fkJ_{\defnu+2,3\defnu+2}$))
and $\barp(j):=1$ ($j\in\fkJ_{\defnu+1,\defN}$ (resp. $\fkJ_{\defnu+1,3\defnu}$,
resp. $\fkJ_{1,\defnu+1}\cup\fkJ_{3\defnu+3,4\defnu+2}$)).
Then we have $\fksl_\barp^{(\deftau)}=\fksl^{(\deftau)}(\defM|\defM)$.
Let $\bard_i:=(-1)^{\barp(i)}$ for $i\in\fkJ_{1,\defN}$.
Define the $3(\hatdefN+1)$ elements $\barh_i$, $\bare_i$, $\barf_i$ ($i\in\bI$)
of $\fksl_\barp^{(\deftau)}$ as follows. \newline\par 
If $\deftau=1$, let 
$\barh_0:=-\bard_1\mate_{11}+\bard_\defN\mate_{\defN\defN}+\checkc$,
$\bare_0:=\bard_N\mate_{\defN1}\otimes t$, $\barf_0:=\mate_{1\defN}\otimes {t^{-1}}$
and $\barh_i:=\bard_i\mate_{ii}-\bard_{i+1}\mate_{i+1i+1}$,
$\bare_i:=\bard_i\mate_{ii+1}$, $\barf_i:=\mate_{i+1i}$ ($i\in\fkJ_{1,\hatdefN}$). \par
If $\deftau=2$, let 
$\barh_0:=-2\bard_1(\mate_{11}-\mate_{\defgamma_\defNp(1)\defgamma_\defNp(1)})+2\checkc$,
$\bare_0:=2\bard_1\mate_{\defgamma_\defNp(1)1}\otimes t$,
$\barf_0:=\mate_{1\defgamma_\defNp(1)}\otimes t^{-1}$,
$\barh_i:=\bard_i\mate_{ii}-\bard_{i+1}\mate_{i+1i+1}
+\bard_{i+1}\mate_{\defgamma_\defNp(i+1)\defgamma_\defNp(i+1)}
-\bard_i\mate_{\defgamma_\defNp(i)\defgamma_\defNp(i)}$,
$\bare_i:=\bard_i(\mate_{ii+1}-(-1)^{\barp(i)\barp(i+1)+\barp(i+1)}\defg^\barp(i)\defg^\barp(i+1)
\mate_{\defgamma_\defNp(i+1)\defgamma_\defNp(i)})$, \newline
$\barf_i:=\mate_{i+1i}-(-1)^{\barp(i+1)\barp(i)+\barp(i)}\defg^\barp(i+1)\defg^\barp(i)
\mate_{\defgamma_\defNp(i)\defgamma_\defNp(i+1)}$ ($i\in\fkJ_{1,\hatdefN-1})$, and \newline
$\barh_\hatdefN:=2\bard_\defN(\mate_{\defN\defN}-\mate_{\defgamma_\defNp(\defN)\defgamma_\defNp(\defN)})$, 
$\bare_\hatdefN:=2\bard_\defN\mate_{\defN\defgamma_\defNp(\defN)}$, and  
$\barf_\hatdefN:=\mate_{\defgamma_\defNp(\defN)\defN}$. \par
If $\deftau=4$, let
$\barh_0:=-\bard_{\deftheta(1)}(\mate_{\deftheta(1)\deftheta(1)}
-\mate_{\deftheta(\defgamma_{\defNp-1}(1))\deftheta(\defgamma_{\defNp-1}(1))})+2\checkc$,
$\bare_0:=-(\mate_{\deftheta(0)\deftheta(1)}-\defg^{\prime,\barp}(\defgamma_{\defNp-1}(1))
\mate_{\deftheta(\defgamma_{\defNp-1}(1))\deftheta(0)})\otimes t$, \newline
$\barf_0:=(\mate_{\deftheta(1)\deftheta(0)}+\defg^{\prime,\barp}(1)
\mate_{\deftheta(0)\deftheta(\defgamma_{\defNp-1}(1))})\otimes t^{-1}$, and 
$\barh_i:=\bard_i\mate_{\deftheta(i)\deftheta(i)}$ \newline $-\bard_{i+1}\mate_{\deftheta(i+1)\deftheta(i+1)}
+\bard_{i+1}\mate_{\defgamma_{\defNp-1}(\deftheta(i+1))\defgamma_{\defNp-1}(\deftheta(i+1))}
-\bard_i\mate_{\defgamma_{\defNp-1}(\deftheta(i))\defgamma_{\defNp-1}(\deftheta(i))}$ 
($i\in\fkJ_{1,\hatdefN-1}$), 
$\barh_\hatdefN:=\bard_\defN\mate_{\deftheta(\defN)\deftheta(\defN)}
-\bard_\defN\mate_{\defgamma_{\defNp-1}(\deftheta(\defN))\defgamma_{\defNp-1}(\deftheta(\defN))}$, \newline
$\bare_i:=\bard_i(\mate_{\deftheta(i)\deftheta(i+1)}
-(-1)^{\barp(\deftheta(i))\barp(\deftheta(i+1))+\barp(\deftheta(i+1))}
\defg^{\prime,\barp}(i)\defg^{\prime,\barp}(i+1)
\mate_{\deftheta(\defgamma_{\defNp-1}(i+1))\deftheta(\defgamma_{\defNp-1}(i))})$, 
$\barf_i:=$ 
$\mate_{\deftheta(i+1)\deftheta(i)}
-(-1)^{\barp(\deftheta(i+1))\barp(\deftheta(i))+\barp(\deftheta(i))}
\defg^{\prime,\barp}(i+1)\defg^{\prime,\barp}(i)
\mate_{\deftheta(\defgamma_{\defNp-1}(i))\deftheta(\defgamma_{\defNp-1}(i+1))}$
($i\in\fkJ_{1,\hatdefN}$).
\newline\newline
Let $\defV$ be an $(\defN+2)$-dimensional $\bR$-linear space 
with a basis $\{\bveps_i|i\in\fkJ_{1,\defN}\}\cup\{\defhatdel,\defprt\}$.
Let $(\,,\,):\defV\times\defV\to\bR$ be the symmetric bilinear map 
defined by $(\bveps_i,\bveps_j):=\delta_{ij}\bard_i$,
$(\bveps_i,\defhatdel):=0$
($i,j\in\fkJ_{1,\defN}$), $(\defhatdel,\defhatdel):=0$,
$(\defprt,\defhatdel):=1$, $(\defprt,\defprt):=0$.
Define $\al_i\in\defV$ ($i\in\bI$) as follows.
Let $\al_0:=\defhatdel-\bveps_1+\bveps_\hatdefN$
(resp. $\defhatdel-2\bveps_1$, resp. $\defhatdel-2\bveps_1$)
if $\deftau=1$ (resp. $2$, resp. $4$).
Let $\al_i:=\bveps_i-\bveps_{i+1}$ for $i\in\fkJ_{1,\defN-1}$.
Let $\al_\defN:=\bveps_N$ (resp. $2\bveps_N$) if $\deftau=2$ (resp. $4$). 
\newline\newline
{\bf{Presentation of $\fksl_\barp^{(\deftau)}$.}} The elements $\checkd$, $\barh_i$, $\bare_i$, $\barf_i$ ($i\in\bI$) generate 
$\fksl_\barp^{(\deftau)}$ as a Lie superalgebra, where
$\checkd$, $\barh_i\in\fksl_\barp^{(\deftau)}(0)$,
and  $\bare_i$, $\barf_i\in\fksl_\barp^{(\deftau)}(\delta_{i\defnu}+\delta_{\deftau 1}\delta_{i0})$. We have the following equations.
\begin{equation}\label{eqn:defrelofsltau}
\begin{array}{l}
\mbox{$[\checkd,\barh_i]=0$, $[\checkd,\bare_i]=\delta_{i0}\bare_i$, $[\checkd,\barf_i]=-\delta_{i0}\barf_i$,} \\
\mbox{$[\barh_i,\barh_j]=0$, $[\barh_i,\bare_j]=(\al_i,\al_j)\bare_j$, $[\barh_i,\barf_j]=-(\al_i,\al_j)\barf_j$,
$[\bare_i,\barf_j]=\delta_{ij}\barh_i$,} \\
\mbox{$[\bare_i,\bare_j]=[\barf_i,\barf_j]=0$ if $(\al_i,\al_j)=0$,} \\
\mbox{$[\bare_i,[\bare_i,\bare_j]]=[\barf_i,[\barf_i,\barf_j]]=0$ if $i\ne j$ and $-2(\al_i,\al_j)=(\al_i,\al_i)\ne 0$,} \\
\mbox{$[\bare_i,[\bare_i,[\bare_i,\bare_j]]]=[\barf_i,[\barf_i,[\barf_i,\barf_j]]]=0$ 
if $i\ne j$ and $-(\al_i,\al_j)=(\al_i,\al_i)\ne 0$,} \\
\mbox{$[[\bare_i,\bare_j],[\bare_i,\bare_k]]=[[\barf_i,\barf_j],[\barf_i,\barf_k]]=0$} \\
\quad\quad\mbox{if $i\ne j\ne k\ne i$, $(\al_i,\al_i)=(\al_j,\al_k)=0$ and $-(\al_i,\al_j)=(\al_i,\al_k)\ne 0$.} 
\end{array}
\end{equation} 
\begin{theorem}\label{theorem:presentation}
{\rm{(\cite[Theorem~4.1.1]{Y99}, \cite{Y01})}}
Assume that $(\deftau,\defM)\ne (1,2)$. Then the equations of \eqref{eqn:defrelofsltau} are the defining relations of $\fksl_\barp^{(\deftau)}$
for the generators $\checkd$, $\barh_i$, $\bare_i$, $\barf_i$ $(i\in\bI)$.
\end{theorem}
Let $\fkntauo$, (resp. $\fkntaup$), (resp. $\fkntaum$) be the sub-Lie superalgebra of $\fksl_\barp^{(\deftau)}$
generated by the elements $\checkd$, $\barh_i$ (resp. $\bare_i$, resp. $\barf_i$) ($i\in\bI$).
Then $\fksl_\barp^{(\deftau)}=\fkntaum\oplus\fkntauo\oplus\fkntaup$ as a $\bK$-linear space,
and the defining relations of $\fkntauo$, $\fkntaup$, $\fkntaum$ for the above generators are those of 
\eqref{eqn:defrelofsltau} 
respectively.  
\newline\newline
{\bf{Dynkin diagrams of $\fksl_\barp^{(\deftau)}=\fksl^{(\deftau)}(\defM|\defM)$.}}
$$
\DiagramAA
$$
$$
\DiagramCC
$$
$$
\DiagramBB
$$

\section{Quantum superalgebras}
\subsection{Multiparameter quantum affine superalgebras $\defUbhm$ and $\defUntbhm$ of $\fksl^{(\deftau)}(M|M)$}
Keep the notation of Subsection~\ref{subsection:sltau}.
Let $\bKt:=\bK\setminus\{0\}$.
Let $\hatdefV:=(\oplus_{i=0}^\hatdefN\bZ\al_i)\oplus\bZ\defprt(\subset\defV)$.
Let $\defq\in\bKt$ be such that $\defq^m\ne 1$ for all $m\in\bN$. 
Let $\bhm:\hatdefV\times\hatdefV\to\bKt$ be a map such that
\begin{equation}\label{eqn:defeqbhm}
\begin{array}{l}
\bhm(a+b,c)=\bhm(a,c)\bhm(b,c),\,\,\bhm(a,b+c)=\bhm(a,b)\bhm(a,c)\quad(a,b,c\in\hatdefV), \\
\bhm(\al_i,\al_i):=\defq^{(\al_i,\al_i)}\quad(i\in\bI,\,(\al_i,\al_i)\ne 0), \\
\bhm(\al_i,\al_i):=-1\quad(i\in\bI,\,(\al_i,\al_i)=0), \\
\bhm(\al_i,\al_j)\bhm(\al_j,\al_i):=\defq^{2(\al_i,\al_j)}
\quad(i,j\in\bI,\,i\ne j), \\
\bhm(\defprt,\defprt):=1,\,\,\bhm(\defprt,\al_i):=\bhm(\al_i,\defprt):=\defq^{\delta_{i0}}.
\end{array}
\end{equation}
Define the map $\bhmop:\hatdefV\times\hatdefV\to\bKt$ by
$\bhmop(a,b):=\bhm(b,a)$ ($a$, $b\in\hatdefV$).
Denote $\hatbhm:\hatdefV\times\hatdefV\to\bKt$ be the map with the same equations as in \eqref{eqn:defeqbhm}
and the equations $\hatbhm(\al_i,\al_j)=\defq^{(\al_i,\al_j)}$ ($i,j\in\bI$,
$i\ne j$).  

Let $\tildefUbhm$ be the associative $\bK$-algebra (with $1$) defined by
the generators $\deftilK_a$, $\deftilL_a$ ($a\in\hatdefV$),
$\deftilE_i$, $\deftilF_i$ ($i\in\bI$), and the relations:
\begin{equation}\label{eqn:bhmfundrl}
\begin{array}{l}
\deftilK_0=\deftilL_0=1,\,
\deftilK_a\deftilK_b=\deftilK_{a+b},\,\deftilL_a\deftilL_b=\deftilL_{a+b},\,\deftilK_a\deftilL_b=\deftilL_b\deftilK_a, \\
\deftilK_a\deftilE_i\deftilK_a^{-1}=\bhm(a,\al_i)\deftilE_i,\,
\deftilK_a\deftilF_i\deftilK_a^{-1}=\bhm(a,-\al_i)\deftilF_i,\\
\deftilL_a\deftilE_i\deftilL_a^{-1}=\bhm(-\al_i,a)\deftilE_i,\,
\deftilL_a\deftilF_i\deftilL_a^{-1}=\bhm(\al_i,a)\deftilF_i, \\
\deftilE_i\deftilF_j-\deftilF_j\deftilE_i=\delta_{ij}(-\deftilK_{\al_i}+\deftilL_{\al_i}) \\
(a,b\in\hatdefV,\,i,j\in\bI).
\end{array}
\end{equation} 
Define the $\bK$-algebra isomorphism $\tilUpsbhmop:\tildefUbhmop\to\tildefUbhm$ by 
$\tilUpsbhmop(\deftilK_a):=\deftilL_a$, \newline 
$\tilUpsbhmop(\deftilL_a):=\deftilK_a$, 
$\tilUpsbhmop(\deftilE_i):=\deftilF_i$, $\tilUpsbhmop(\deftilF_i):=\deftilE_i$
($a\in\hatdefV$, $i\in\bI$).

We see that $\tildefUbhm$ is the Hopf algebra $(\tildefUbhm,\defDelta,\defS,\defepsilon)$
with 
$\defDelta(\deftilK_a)=\deftilK_a\otimes\deftilK_a$, $\defDelta(\deftilL_a)=\deftilL_a\otimes\deftilL_a$, 
$\defDelta(\deftilE_i)=\deftilE_i\otimes 1 + \deftilK_{\al_i}\otimes\deftilE_i$, 
$\defDelta(\deftilF_i)=\deftilF_i\otimes\deftilL_{\al_i} + 1\otimes\deftilF_i$, 
$\defS(\deftilK_a)=\deftilK_a^{-1}$, $\defS(\deftilL_a)=\deftilL_a^{-1}$, 
$\defS(\deftilE_i)=-\deftilK_{\al_i}^{-1}\deftilE_i$, $\defS(\deftilF_i)=-\deftilF_i\deftilL_{\al_i}^{-1}$, 
$\defepsilon(\deftilK_a)=\defepsilon(\deftilL_a)=1$,
$\defepsilon(\deftilE_i)=\defepsilon(\deftilF_i)=0$
$(a,b\in\hatdefV,\,i,j\in\bI)$.

Let $\tildefUobhm$ (resp. $\tildefUpbhm$, resp. $\tildefUmbhm$)
be the $\bK$-subalgebra of $\tildefUbhm$
generated by the elements $\deftilK_a$, $\deftilL_a$ ($a\in\hatdefV$)
(resp. $1$, $\deftilE_i$ ($i\in\bI$), resp. $1$, $\deftilF_i$ ($i\in\bI$)).
Then the elements $\deftilK_a\deftilL_b$ ($a$, $b\in\hatdefV$) form a $\bK$-basis of $\tildefUobhm$,
and $\tildefUpbhm$ (resp. $\tildefUmbhm$) is the free $\bK$-algebra in $\deftilE_i$'s
(resp. $\deftilF_i$'s). Moreover we have the $\bK$-linear isomorphism
$\tildefvarthetabhm:\tildefUmbhm\otimes\tildefUobhm\otimes\tildefUpbhm\to\tildefUbhm$
defined by $\tildefvarthetabhm(x\otimes y\otimes z):=xyz$.
Let $\hatdefVp:=\oplus_{i\in\bI}\bZ\al_i$ and
$\hatdefVprp:=\oplus_{i\in\bI}\bZgeqo\al_i$.
Define the $\bK$-subspaces $(\tildefUbhm)_\lambda$
($\lambda\in\hatdefVp$) by $\tildefUbhm=\oplus_{\lambda\in\hatdefVp}(\tildefUbhm)_\lambda$,
$\tildefUobhm\subset(\tildefUbhm)_0$,
$\deftilE_i\in(\tildefUbhm)_{\al_i}$ and
$\deftilF_i\in(\tildefUbhm)_{-\al_i}$. 
Let $(\tildefUpbhm)_\lambda:=\tildefUpbhm\cap(\tildefUbhm)_\lambda$
and $(\tildefUmbhm)_\lambda:=\tildefUmbhm\cap(\tildefUbhm)_\lambda$
($\lambda\in\hatdefVp$). Then
$\tildefUpbhm=\oplus_{\lambda\in\hatdefVprp}(\tildefUpbhm)_\lambda$,
$\tildefUmbhm=\oplus_{\lambda\in\hatdefVprp}(\tildefUmbhm)_{-\lambda}$
and $(\tildefUpbhm)_0=(\tildefUmbhm)_0=\bK\cdot 1$.
For each $\lambda\in\hatdefVp$, define 
the $\bK$-subspace $(\tilIpbhm)_\lambda$ of $(\tildefUpbhm)_\lambda$ as follows.
Let $(\tilIpbhm)_\lambda:=\{0\}$ if $\lambda=0$ or $\lambda\notin\hatdefVprp$.
If $\lambda\in\hatdefVprp\setminus\{0\}$,
let $(\tilIpbhm)_\lambda$ be formed by all the elements
$X\in(\tildefUpbhm)_\lambda$ such that
\begin{equation}
\forall i\in\bI,\quad X\deftilF_i-\deftilF_iX\in(\tilIpbhm)_{\lambda-\al_i}\deftilK_{\al_i}+(\tilIpbhm)_{\lambda-\al_i}\deftilL_{\al_i}.
\end{equation}
Let $\tilIpbhm:=\oplus_{\lambda\in\hatdefVprp}(\tilIpbhm)_\lambda$,
and $\tilImbhm:=\tilUpsbhmop(\tilIpbhmop)$.
Then $\tilIpbhm$ (resp. $\tilImbhm$)  is an ideal of  $\tildefUpbhm$ (resp. $\tildefUmbhm$).
Let $\tilJbhm:=\rmSpan_\bK(\tildefUmbhm\tildefUobhm\tilIpbhm)
+\rmSpan_\bK(\tilImbhm\tildefUobhm\tildefUpbhm)$.
Then $\tilJbhm$ is a Hopf-ideal of  $\tildefUbhm$
(see \cite[\S~4.1]{AYY15} for example). Let $\defUbhm$ denote
the Hopf $\bK$-algebra defined by $\defUbhm:=\tildefUbhm/\tilJbhm$.

Let $\tilpibhm:\tildefUbhm\to\defUbhm$ be the canonical map.
Let $\defK_a:=\tilpibhm(\deftilK_a)$, $\defL_a:=\tilpibhm(\deftilL_a)$
($a\in\hatdefVp$), 
and $\defE_i:=\tilpibhm(\deftilE_i)$, $\defF_i:=\tilpibhm(\deftilF_i)$
($i\in\bI$). 
Let $\defUobhm:=\tilpibhm(\tildefUobhm)$,
$\defUpbhm:=\tilpibhm(\tildefUpbhm)$ and $\defUmbhm:=\tilpibhm(\tildefUmbhm)$.

\begin{remark} \label{remark:NoteYandEG}
Let $X=(X,\Delta_X, S_X, \epsilon_X)$ be a $\bK$-Hopf algebra.
Assume that
$X=X(0)\oplus X(1)$ as a $\bK$-linear space.
For $t\in\bZ\setminus\fkJ_{0,1}$ with
$t=2t_1+t_2$ for some $t_1\in\bZ$ and $t_2\in\fkJ_{0,1}$,
let $X(t):=X(t_2)$.
Assume that
$X(s)X(t)\subset X(s+t)$ ($s$, $t\in\bZ$)
and $\Delta_X(X(t))\subset X(0)\otimes X(t)+X(1)\otimes X(t-1)$
($t\in\bZ$).
Then we define a $\bK$-Hopf algebra
$X^\defsigma=(X^\defsigma,\Delta_{X^\defsigma}, S_{X^\defsigma}, \epsilon_{X^\defsigma})$
with the following conditions (1) and (2). (1) $X^\defsigma$ contains $X$
as a $\bK$-Hopf subalgebra.
(2) There exists 
$\defsigma\in X^\defsigma$ so that
$\defsigma^2=1$, $\Delta_{X^\defsigma}(\defsigma)=\defsigma\otimes\defsigma$,
$S_{X^\defsigma}(\defsigma)=\defsigma$, $\epsilon_{X^\defsigma}(\defsigma)=1$,
$\defsigma x \defsigma
=(-1)^t x$ ($t\in\fkJ_{0,1}$, $x\in X(t)$)
and $X^\defsigma=X\otimes\defsigma X$ as a $\bK$-linear space.
Define the map $\defwp:\bI\to\fkJ_{0,1}$ by
$\defwp(i):=1-\delta_{1,(\al_i,\al_i)}$.
Define the group homomorphism $\hatdefwp:\hatdefVp\to\{-1,1\}(\subset\bZ)$
by  $\hatdefwp(\al_i):=(-1)^{\hatdefwp(i)}$ ($i\in\bI$).
Let $\defUbhm^\defsigma$ be the one defined for
the direct sum $\defUbhm=\defUbhm(0)\oplus\defUbhm(1)$
with $\defUbhm(t):=\oplus_{\lambda\in\hatdefwp^{-1}(\{(-1)^t\})}(\defUbhm)_\lambda$
($t\in\fkJ_{0,1}$).
Let $\defUbhm^{\defsigma,\prime}$ be the $\bK$-subalgebra
of $\defUbhm^\defsigma$ generated by
$\defK_{\al_i}^{\pm 1}$, $\defL_{\al_i}^{\pm 1}\defsigma^{\defwp(i)}$,
$\defE_i$, $\defF_i\defsigma^{\defwp(i)}$ ($i\in\bI$)
and $\defK_\defprt^{\pm 1}$, $\defL_\defprt^{\pm 1}$.
Assume that $\bhm(\al_i,\al_j)=\defq^{(\al_i,\al_j)}$ ($i$, $j\in\bI$).
Let $Y$ be the ideal (as a $\bK$-algebra) of $\defUbhm^{\defsigma,\prime}$ 
generated by $\defK_{\al_i}\defL_{\al_i}\defsigma^{\defwp(i)}-1$
($i\in\bI$) and $\defK_\defprt\defL_\defprt-1$.
Let $\defUbhm^\sharp:=\defUbhm^{\defsigma,\prime}/Y$ as the quotient $\bK$-algebra.
Recall the notation  $Z^{(\deftau)}_k$ from Introduction. From 
Subsection~\ref{subsection:defAtaumn}, we see
\begin{equation*}
\rmA^{(\deftau)}(\defM-1,\defM-1)=\fksl^{(\deftau)}(\defM|\defM))/(\oplus_{k=-\infty}^\infty\bC Z^{(\deftau)}_k)
\quad(\deftau\in\{1,2,4\}).
\end{equation*}
If $\deftau\in\{2,4\}$, 
$\defU_q(\rmA^{(\deftau)}(\defM-1,\defM-1))$
is isomorphic to $\defUbhm^\sharp$ as a $\bK$-algebra.
Let $\gamma:=\sum_{t=1}^\defnu t(\al_t+\al_{\defN-t})\in\hatdefVp$. Then $\bZ\gamma=\{\lambda\in\hatdefVp|\forall \mu\in\hatdefVp,
(\lambda,\mu)=0\}$.
 If $\deftau=1$, $\defUbhm^\sharp/(\defK_\gamma-1)\defUbhm^\sharp$
(resp. $\defUbhm^\sharp$)
is isomorphic to $\defU_q(\rmA^{(1)}(\defM-1,\defM-1))$
(resp. $\defU_q(\fksl^{(1)}(\defM|\defM))/(\oplus_{k\ne 0}\bC Z^{(1)}_k)$).
See also \cite[\S 8]{Y94}, \cite{Y01}, \cite{EG09}, \cite[Remark~7.11]{AYY15}.
\end{remark}

\begin{theorem}\label{theorem:weightA}
{\rm{(\cite{EG09}, see also \cite[Theorem~8.4.3]{Y94}, \cite{Y01})}}
Assume that $\defq$ to be transcendental over $\bQ$.
Assume that $\bhm(\al_i,\al_j)=\defq^{(\al_i,\al_j)}$ for all $i$, $j\in\bI$ with $i\ne j$.
Then for each $\lambda\in\hatdefVprp$,
$\dim(\defUpbhm)_\lambda$ equals the dimension of the weight space $\cU(\brvpi^{(1)}(\fkntaup))_\lambda$ of $\lambda$ of
the universal enveloping superalgebra $\cU(\brvpi^{(1)}(\fkntaup))$ of $\brvpi^{(1)}(\fkntaup)(\subset\rmA^{(\deftau)}_\barp
=\rmA^{(\deftau)}(\defM-1,\defM-1))$.
\end{theorem}

For $\lambda$, $\mu\in\hatdefVp$
and ${\tilde X}\in(\tildefUbhm)_\lambda$, ${\tilde Y}\in(\tildefUbhm)_\mu$, 
let $\wbrl {\tilde X},{\tilde Y}\wbrr:={\tilde X}{\tilde Y}-{\frac 1 {\bhm(\mu,\lambda)}}{\tilde Y}{\tilde X}
(\in(\tildefUbhm)_{\lambda+\mu})$.
Let $\tilIntpbhm$ be the ideal of the $\bK$-algebra $\tildefUpbhm$ generated by the following elements:
\begin{equation*}
\begin{array}{l}
\mbox{$\wbrl\deftilE_i,\deftilE_j\wbrr$  ($(\al_i,\al_j)=0$),}\,\,
\mbox{$\wbrl\deftilE_i,\wbrl\deftilE_i,\deftilE_j\wbrr\wbrr$ ($i\ne j$, $-2(\al_i,\al_j)=(\al_i,\al_i)\ne 0$),} \\
\mbox{$\wbrl\deftilE_i,\wbrl\deftilE_i,\wbrl\deftilE_i,\deftilE_j\wbrr\wbrr\wbrr$ 
($i\ne j$, $-(\al_i,\al_j)=(\al_i,\al_i)\ne 0$),} \\
\mbox{$\wbrl\wbrl\deftilE_i,\deftilE_j\wbrr,\wbrl\deftilE_i,\deftilE_k\wbrr\wbrr$} \\
\quad\quad\mbox{($i\ne j\ne k\ne i$, $(\al_i,\al_i)=(\al_j,\al_k)=0$, $-(\al_i,\al_j)=(\al_i,\al_k)\ne 0$).} 
\end{array}
\end{equation*} Then we can see $\tilIntpbhm\subset\tilIpbhm$.
Let $\tilIntmbhm:=\tilUpsbhmop(\tilIntpbhmop)$. Then $\tilIntmbhm\subset\tilImbhm$.
Let $\tilJntbhm:=\rmSpan_\bK(\tildefUmbhm\tildefUobhm\tilIntpbhm)
+\rmSpan_\bK(\tilIntmbhm\tildefUobhm\tildefUpbhm)$.
Then $\tilJntbhm$ is a Hopf-ideal of  $\tildefUbhm$
(see \cite[\S~4.1]{AYY15} for example). Let $\defUntbhm$ denote
the Hopf $\bK$-algebra defined by $\defUntbhm:=\tildefUbhm/\tilJntbhm$.
Let $\tilpintbhm:\tildefUbhm\to\defUntbhm$ be the canonical map.
Let $\defntK_a:=\tilpintbhm(\deftilK_a)$, $\defntL_a:=\tilpintbhm(\deftilL_a)$
($a\in\hatdefVp$), 
and $\defntE_i:=\tilpintbhm(\deftilE_i)$, $\defntF_i:=\tilpintbhm(\deftilF_i)$
($i\in\bI$). 
Then we have the Hopf $\bK$-algebra epimorphism
$\pintbhm:\defUntbhm\to\defUbhm$ such that
$\pintbhm(\defntK_a)=\defK_a$, $\pintbhm(\defntL_a)=\defL_a$
($a\in\hatdefVp$), 
and $\pintbhm(\defntE_i)=\defE_i$, $\pintbhm(\defntF_i)=\defF_i$
($i\in\bI$). 
Let $\defUntobhm:=\tilpintbhm(\tildefUobhm)$,
$\defUntpbhm:=\tilpintbhm(\tildefUpbhm)$ and
$\defUntmbhm:=\tilpintbhm(\tildefUmbhm)$.
Then we have the $\bK$-linear isomorphism
$\defUntmbhm\otimes\defUntobhm\otimes\defUntpbhm\to\defUntbhm$
($y\otimes z\otimes x\mapsto xzy$).

\subsection{Vector representation of $\defUntbhm$}
Keep the notation of Subsection~\ref{subsection:sltau}.
We denote by ${\mathrm{M}}_\defNp(\bK)$ the $\bK$-algebra
such that it is $\fkgl_\barp$ as a $\bK$-linear space and 
its multiplication is defined by $\mate_{ij}\mate_{kl}=\delta_{jl}\mate_{ij}$
($i$, $j$, $k$, $l\in\fkJ_{1,\defNp}$). Let $\{\mate_i|i\in\fkJ_{1,\defNp}\}$
be the natural $\bK$-basis of $\bK^\defNp$.
We regard ${\mathrm{M}}_\defNp(\bK)\otimes_\bK\bK[t,t^{-1}]$
as the $\bK$-subalgebra of ${\mathrm{End}}_\bK(\bK^\defNp\otimes_\bK\bK[t,t^{-1}])$ by
$(\mate_{ij}\otimes t^l)(\mate_k\otimes t^m)=\delta_{jk}\mate_i\otimes t^{l+m}$
($i$, $j$, $k\in\fkJ_{1,\defNp}$, $l$, $m\in\bZ$).
\begin{theorem}  
Let $z_i$, $u_i\in\bKt$ {\rm{(}}$i\in\bI${\rm{)}}. 
\begin{equation}\label{eqn:ass}
\begin{array}{l}
\mbox{Assume that $\bhm(\al_i,\defhatdel)=1$
{\rm{(}}resp. $\bhm(\al_i,2\defhatdel)=1${\rm{)}} for all $i\in\bI$} \\
\mbox{if $\deftau\in\{1,2\}$ {\rm{(}}resp. $\deftau=4${\rm{)}}.}
\end{array}
\end{equation}
Then there exists a $\bK$-algebra homomorphism
$\defPsi:\defUntbhm\to{\mathrm{End}}_\bK(\bK^\defNp\otimes_\bK\bK[t,t^{-1}])$
as follows. Let $x_{ij}:=\bhm(\al_i,\al_j)$ {\rm{(}}$i$, $j\in\bI${\rm{)}}.
\newline\newline
$\defPsi(\defK_\defprt)(\mate_k\otimes t^m)=\defq^m\mate_k\otimes t^m$,
$\defPsi(\defL_\defprt)(\mate_k\otimes t^m)=\defq^{-m}\mate_k\otimes t^m$
{\rm{(}}$k\in\fkJ_{1,\defNp}$, $m\in\bZ${\rm{)}}. \newline
For $\deftau=1$, $\defPsi(\defntK_{\al_i})=z_i\sum_{t=1}^\defN(\prod_{s=0}^{t-1}x_{is}^{-1})\mate_{tt}$
{\rm{(}}$i\in\bI${\rm{)}}, \newline
$\defPsi(\defntL_{\al_i})=z_i\defq^{-2(\bveps_\defN,\al_i)}\sum_{t=1}^\defN(\prod_{s=0}^{t-1}x_{si})\mate_{tt}$
{\rm{(}}$i\in\bI${\rm{)}}, 
$\defPsi(\defntE_0)=z_0u_0\bard_1\defq^{\bard_1}(\defq-\defq^{-1})\mate_{\defN 1}\otimes t$,
$\defPsi(\defntE_i)=-z_iu_i\bard_i\defq^{-\bard_i}(\prod_{s=0}^{i-1}x_{is}^{-1})(\defq-\defq^{-1})\mate_{ii+1}$
{\rm{(}}$i\in\bI\setminus\{0\}${\rm{)}},
$\defPsi(\defntF_0)=u_0^{-1}\mate_{1\defN}\otimes t^{-1}$,
$\defPsi(\defntF_i)=u_i^{-1}\mate_{i+1i}$ {\rm{(}}$i\in\bI\setminus\{0\}${\rm{)}}.
\newline
For $\deftau=2$, \newline
$\defPsi(\defntK_{\al_i})=z_i((\sum_{t=1}^\defN(\prod_{s=0}^{t-1}x_{is}^{-1})\mate_{tt})+
(\sum_{t=1}^\defN(\prod_{s=0}^{t-1}x_{is})\mate_{\defgamma_\defNp(t)\defgamma_\defNp(t)}))$
{\rm{(}}$i\in\bI${\rm{)}},  
$\defPsi(\defntL_{\al_i})=z_i\defq^{2(\bveps_1,\al_i)}((\sum_{t=1}^\defN(\prod_{s=0}^{t-1}x_{si})\mate_{tt})+
(\sum_{t=1}^\defN(\prod_{s=0}^{t-1}x_{si}^{-1})\mate_{\defgamma_\defNp(t)\defgamma_\defNp(t)}))$
{\rm{(}}$i\in\bI${\rm{)}}, 
$\defPsi(\defntE_0)=-z_0u_0\bard_1\defq^{-2\bard_1}(\defq^2-\defq^{-2})\mate_{\defgamma_\defNp(1)1}\otimes t$,
$\defPsi(\defntE_i)=$ \newline $-z_iu_i\bard_i\defq^{-\bard_i}(\prod_{s=0}^{i-1}x_{is}^{-1})(\defq-\defq^{-1})
(\mate_{ii+1}-\defq^{4\bard_1+2\bard_i}(\prod_{s=1}^{i-1}x_{is}^2)\mate_{\defgamma_\defNp(i+1)\defgamma_\defNp(i)})$
{\rm{(}}$i\in\bI\setminus\{0,\defN\}${\rm{)}}, 
$\defPsi(\defntE_\defN)=-z_\defN u_\defN\bard_\defN(\prod_{s=1}^{\defN-1}x_{\defN s}^{-1})
\defq^{-2\bard_\defN}(\defq^2-\defq^{-2})\mate_{\defN\defgamma_\defNp(\defN)}$, \newline
$\defPsi(\defntF_0)=u_0^{-1}\mate_{1\defgamma_\defNp(1)}\otimes t^{-1}$,
$\defPsi(\defntF_i)=u_i^{-1}(\mate_{i+1i}-\mate_{\defgamma_\defNp(i)\defgamma_\defNp(i+1)})$
{\rm{(}}$i\in\bI\setminus\{0,\defN\}${\rm{)}},
$\defPsi(\defntF_\defN)=u_\defN^{-1}\mate_{\defgamma_\defNp(\defN)\defN}$.
\newline
For $\deftau=4$, \newline
$\defPsi(\defntK_{\al_i})=z_i((\sum_{t=0}^{\defN+1}(\prod_{s=0}^{t-1}x_{is}^{-1})\mate_{\deftheta(t)\deftheta(t)})$
\newline
$+
(\sum_{t=1}^\defN(\prod_{s=0}^{t-1}x_{is})\mate_{\deftheta(\defgamma_{\defNp-1}(t))\deftheta(\defgamma_{\defNp-1}(t))}))$
{\rm{(}}$i\in\bI${\rm{)}}, \newline 
$\defPsi(\defntL_{\al_i})=z_i((\sum_{t=0}^{\defN+1}(\prod_{s=0}^{t-1}x_{si})\mate_{\deftheta(t)\deftheta(t)})
$ \newline
$+
(\sum_{t=1}^\defN(\prod_{s=0}^{t-1}x_{si}^{-1})\mate_{\deftheta(\defgamma_{\defNp-1}(t))\deftheta(\defgamma_{\defNp-1}(t))}))$
{\rm{(}}$i\in\bI${\rm{)}}, \newline
$\defPsi(\defntE_0)=-z_0u_0\bard_1(\defq-\defq^{-1})(\mate_{\deftheta(0)\deftheta(1)}
-\mate_{\deftheta(\defgamma_{\defNp-1}(1))\deftheta(0)})\otimes t$, \newline
$\defPsi(\defntE_i)=-z_iu_i\bard_i\defq^{-\bard_i}(\prod_{s=0}^{i-1}x_{is}^{-1})(\defq-\defq^{-1})
(\mate_{\deftheta(i)\deftheta(i+1)}$ \newline $-\defq^{2\bard_i}(\prod_{s=0}^{i-1}x_{is}^2)\mate_{\defgamma_{\defNp-1}(\deftheta(i+1))\defgamma_{\defNp-1}(\deftheta(i))})$
{\rm{(}}$i\in\bI\setminus\{0\}${\rm{)}}, \newline
$\defPsi(\defntF_0)=u_0^{-1}(\mate_{\deftheta(1)\deftheta(0)}
-\mate_{\deftheta(0)\deftheta(\defgamma_{\defNp-1}(1))})\otimes t^{-1}$, \newline
$\defPsi(\defntF_i)=u_i^{-1}
(\mate_{\deftheta(i+1)\deftheta(i)}-\mate_{\defgamma_{\defNp-1}(\deftheta(i))\defgamma_{\defNp-1}(\deftheta(i+1))})$
{\rm{(}}$i\in\bI\setminus\{0\}${\rm{)}}. 
\end{theorem}
{\it{Proof.}} This can be proved directly. \hfill $\Box$

\subsection{Main result---Lowest positive central element of $\defUntbhm$}

Recall $\defnu:={\frac \defN 2}$.
Let $\defs:=1$ (resp. $\defs:=2$) if $\deftau\in\{1,2\}$
 (resp. $\deftau:=4$).
Let $\defntB_{\defnu-1}:=\defntE_\defnu\in(\defUntpbhm)_{\bveps_\defnu-\bveps_{\defnu+1}}$.
For $i\in\fkJ_{0,\defnu-2}$, let
$\defntB_i:=\wbrl\defntE_{\defN-i-1},\wbrl\defntB_{i+1},\defntE_{i+1}\wbrr\wbrr$ 
$\in(\defUntpbhm)_{\bveps_{i+1}-\bveps_{\defN-i}}$.
Let $\defntB_{-1}:=\wbrl\defntE_\defN,\wbrl\defntB_0,\defntE_0\wbrr\wbrr$
if $\deftau\in\{2,4\}$.
Define 
$\defntA_0\in(\defUntpbhm)_{\defs\defhatdel-\bveps_1+\bveps_\defN}$ to be
$\defntE_0$ (resp. $\defntB_{-1}$,
resp. $\wbrl\defntE_\defN,\wbrl\defntB_{-1},\defntE_0\wbrr\wbrr$)
if $\deftau=1$ (resp. $\deftau=2$, resp. $\deftau=4$).
For $i\in\fkJ_{1,\defnu-1}$, let $\defntA_i:=\wbrl\defntE_{\defN-i},\wbrl\defntA_{i-1},\defntE_i\wbrr\wbrr
\in(\defUntpbhm)_{\defs\defhatdel-\bveps_{i+1}+\bveps_{\defN-i}}$.
Let $\defa_0:=1\in\bKt$.
For $i\in\fkJ_{1,\defnu-2}$, let \newline
$\defa_i:=
a_{i-1}\cdot\bhm(\defs\defhatdel,\al_{\defN-i})\bhm(-\al_i,\sum_{t={i+1}}^{\defN-i-1}\al_t)\bhm(\sum_{t=i}^{\defN-i-1}\al_t,-\al_{\defN-i})$. \newline
Let $\defa_{-1}:=\bhm(\al_\defN,\al_0)\bhm(\al_0,\defhatdel)\bhm(\al_\defN,\defhatdel)$
if $\deftau=4$.
Then $\defntZ\in(\defUntpbhm)_{\defs\defhatdel}$ by
\begin{equation}\label{eqn:defZnatural}
\defntZ:=
\left\{
\begin{array}{ll}
\sum_{i=0}^{\defnu-1}\defa_i\wbrl\defntA_i,\defntB_i\wbrr
&\quad\mbox{if $\deftau\in\{1,2\}$}, \\
\defa_{-1}(\defntB_{-1})^2+\sum_{i=0}^{\defnu-1}\defa_i\wbrl\defntA_i,\defntB_i\wbrr
&\quad\mbox{if $\deftau=4$}.
\end{array}\right.
\end{equation}
Then we have our main theorem:
\begin{theorem}\label{theorem:main}
{\rm{(1)}} It follows that
\begin{equation}\label{eqn:profdefZnatural}
\forall i\in\bI,\,\,\defntZ\defntF_i-\defntF_i\defntZ=\wbrl\defntZ,\defntE_i\wbrr=0.
\end{equation}
\newline
{\rm{(2)}} Under the assumption \eqref{eqn:ass}, $\defPsi(\defntZ)=b(\brvee\otimes t^\defs)$
for some $b\in\bKt$,
where $\brvee$ denotes the unit matrix $\sum_{i=1}^\defNp\mate_{ii}$ as above.
\end{theorem}
{\it{Proof.}} The claims can be obtained directly by using calculation formulas similar to those
given in \cite{Y94} and \cite{Y99}. 

Calculations needed for the proof of Theorem~\ref{theorem:main} are almost the same as those of
\cite[\S6]{Y94}. Here we explain it by using an example. Assume $\deftau=2$.
From Subsection~\ref{subsection:sltau}, we have assumed $\defnu\geq 2$ and we have defined $\defN=2\defnu$.
As similar equations to \cite[(4.4.1), (4.4.2)]{Y94},
for $\lambda_t\in\hatdefV$
and $X_t\in(\defUntbhm)_{\lambda_t}$, ($t\in\fkJ_{1,3}$), we have
\eqref{eqn:appone}-\eqref{eqn:appthree} below.
\begin{equation}\label{eqn:appone}
\begin{array}{l}
\wbrl \wbrl X_1,X_2\wbrr,X_3 \wbrr-\wbrl X_1,\wbrl X_2,X_3\wbrr\wbrr \\
\quad = -{\frac 1 {\bhm(\lambda_2,\lambda_1)}}X_2\wbrl X_1,X_3\wbrr
+{\frac 1 {\bhm(\lambda_3,\lambda_2)}}\wbrl X_1,X_3\wbrr X_2.
\end{array}
\end{equation} 
\begin{equation}\label{eqn:apptwo}
\wbrl X_2,X_1\wbrr=-\bhm(\lambda_2,\lambda_1)\wbrl X_1,X_2\wbrr\quad
\mbox{if $\bhm(\lambda_2,\lambda_1)\bhm(\lambda_1,\lambda_2)=1$}.
\end{equation}
\begin{equation}\label{eqn:appthree}
\defS(\wbrl X_1,X_2\wbrr)=\bhm(\lambda_1,\lambda_2)\defntK_{-\lambda_1-\lambda_2}\wbrl X^\prime_2,X^\prime_1\wbrr\quad
\mbox{where $X^\prime_t:=\defntK_{\lambda_t}X_t$}.
\end{equation}

For $i,j\in\fkJ_{1,\defN}$ with $i<j$, define the two elements $\defntE_{\bveps_i-\bveps_j}$,
$\defntE_{\bveps_i+\bveps_j}$ of $\defUntpbhm$ by
$\defntE_{\bveps_i-\bveps_{i+1}}:=\defntE_i$ (if $i\in\fkJ_{1,\defN-1}$), 
$\defntE_{\bveps_i-\bveps_j}:=\wbrl\defntE_{\bveps_i-\bveps_{j-1}},\defntE_{j-1}\wbrr$
 (if $i\in\fkJ_{1,\defN-1}$ and $j\in\fkJ_{i+2,\defN}$), 
$\defntE_{\bveps_i+\bveps_N}:=\wbrl\defntE_{\bveps_i-\bveps_N},\defntE_N\wbrr$
 (if $i\in\fkJ_{1,\defN-1}$),
and  $\defntE_{\bveps_i+\bveps_j}:=\wbrl\defntE_{\bveps_i+\bveps_{j+1}},\defntE_j\wbrr$
 (if $i\in\fkJ_{1,\defN-2}$ and $j\in\fkJ_{i+1,\defN-1}$).

Using \eqref{eqn:appone}, \eqref{eqn:apptwo} and equations similar to those of \cite[Lemma~6.1.1]{Y94}, we can see
\begin{equation}\label{eqn:appfour}
\begin{array}{l}
\wbrl \defntE_i, \defntE_{\bveps_i\pm \bveps_j}\wbrr=0\,(j\geq i+2),\quad
\wbrl \defntE_{\bveps_i- \bveps_j}, \defntE_{j-1} \wbrr=0\,(j\geq i+2), \\
\wbrl \defntE_{\bveps_i+\bveps_j}, \defntE_j \wbrr=0,\quad\wbrl \defntE_{\bveps_i\pm \bveps_j}, \defntE_k \wbrr=0
\,(k\notin\{i-1, i, j-1,j\}).
\end{array}
\end{equation}
By \eqref{eqn:appone}-\eqref{eqn:appfour}, we have
\begin{equation}\label{eqn:appfive}
\begin{array}{l}
\wbrl\defntB_i,\defntE_j\wbrr=0\,(j\notin\{i,\defN-i-1,\defN-i\}),\,
\wbrl\defntE_{\defN-i-1},\defntB_i\wbrr=0,
\\
\wbrl\defntA_i,\defntE_j\wbrr=0\,(j\notin\{i+1,\defN-i-1,\defN-i\}),\,
\wbrl\defntE_{\defN-i},\defntA_i\wbrr=0
\end{array}
\end{equation} for $i\in\fkJ_{0,\defnu-1}$. 
Using these equations, we have $\wbrl\defntZ,\defntE_i\wbrr=0$ for \eqref{eqn:profdefZnatural}.

By an argument similar to that for \eqref{eqn:appfive} and an argument similar to that of \cite[\S6.9-6.11]{Y99}, we see that for `most' $k$ and $i$, lettimg $Y$ be $\defntA_i$ or $\defntB_i$, 
we have $\wbrl Y,\defntF_k\defntK_{-\al_k}\wbrr=0$, which implies $Y\defntF_k-\defntF_kY=0$.
Using this fact, we have $\defntZ\defntF_i-\defntF_i\defntZ=0$ for \eqref{eqn:profdefZnatural}.
\hfill $\Box$

\begin{theorem} Assume that $(\deftau,\defM)\ne (1,2)$
(see Subsection~\ref{subsection:sltau} for notation).
Assume $\defq$ to be transcendental over $\bQ$.  
Assume that for $i$, $j\in\bI$ with $i<j$,
there exists $n_{ij}\in\bZ$ with $\bhm(\al_i,\al_j)=\defq^{n_{ij}}$.
Assume that the condition \eqref{eqn:ass} is fulfilled.
\newline
{\rm{(1)}} For each $\lambda\in\hatdefVprp$,
$\dim(\defUntpbhm)_\lambda$ equals the dimension of the weight space 
$\cU(\fkntaup)_\lambda$ of $\lambda$ of
the universal enveloping superalgebra $\cU(\fkntaup)$
of $\fkntaup(\subset\fksl^{(\deftau)}_\barp
=\fksl^{(\deftau)}(\defM|\defM))$.
{\rm{(}}Here we let $(\defUntpbhm)_\lambda:=\tilpintbhm((\tildefUpbhm)_\lambda)$.{\rm{)}}
\newline
{\rm{(2)}}
The ideal {\rm{(}}as a $\bK$-algebra{\rm{)}} of $\defUbhm$ generated by $\defntZ$
is a Hopf ideal.
{\rm{(}}We have `new' Hopf algebras $\defUntbhm/\defntZ\defUntbhm$
and $\defUntbhm/(\defntZ\defUntbhm+(\defntZ)^\prime\defUntbhm)$,
where $(\defntZ)^\prime$ is the element of $\defUntmbhm$
defined in the same way as that for $\defntZ$.{\rm{)}}
\end{theorem}
{\it{Proof.}}  (1) We obtain the natural representation of $\fksl^{(\deftau)}_\barp$
from $\defPsi$ by taking specialization $\defq\to 1$ and taking conjugation with a diagonal matrix
with diagonal components are $1$ or $-1$.
Incidentally the coproduct (as a Hopf superalgebra) of the universal enveloping superalgebra of Lie superalgebra 
$\fksl^{(\deftau)}_\barp$
is obtained from the coproduct of $\defUntbhm$. 
Then we can prove the claim by a standard argument. 

(2) By the same argument as that for \cite[Lemma~6.6.1]{Y99},
we see that $\dim(\defUpbhm)_\mu=\dim(\defUntpbhm)_\mu$
for $\mu\in\hatdefVprp$ with $\defs\defhatdel-\mu\in\hatdefVprp\setminus\{0\}$.
Then the claim easily follows from 
this fact. 
\hfill $\Box$

\begin{remark}\label{remark:DefUqsltau}
Recall Remark~\ref{remark:NoteYandEG}.  As in Remark~\ref{remark:NoteYandEG}, assume that $\bhm(\al_i,\al_j)=\defq^{(\al_i,\al_j)}$ ($i$, $j\in\bI$).
We define $\defUbhm^{\natural,\sharp}$ 
(resp. $\defUbhm^{\natural,\sharp,\prime}$,
resp. $\defUbhm^{\natural,\sharp,\prime\prime}$) from $\defUntbhm$
(resp. $\defUntbhm/\defntZ\defUntbhm$,
resp. $\defUntbhm/(\defntZ\defUntbhm+(\defntZ)^\prime\defUntbhm)$) 
in the same way as that for $\defUbhm^\sharp$ defined
from $\defUbhm$.
Assume that $(\deftau,\defM)\ne (1,2)$
(see Subsection~\ref{subsection:sltau} for notation).
Then $\defUbhm^{\natural,\sharp}$ 
(resp. $\defUbhm^{\natural,\sharp,\prime}$,
resp. $\defUbhm^{\natural,\sharp,\prime\prime}$) is isomorphic to $U_q(\fksl^{(\deftau)}(M|M))$ 
(resp. $U_q(\fksl^{(\deftau)}(M|M)/Z^{(\deftau)}_1)$,
resp. $U_q(\fksl^{(\deftau)}(M|M)/(Z^{(\deftau)}_{-\delta_{1\deftau}}+Z^{(\deftau)}_1))$)
as a $\bK$-algebra.
If $(\deftau,\defM)\ne (1,2)$, we have a natural epimorphism from 
$\defUbhm^{\natural,\sharp}$ to $U_q(\fksl^{(\deftau)}(M|M))$,
but its kernel must be very big (see also Remark~\ref{remark:smallUqsltau}). 
\end{remark}

\begin{remark}\label{remark:smallUqsltau}
Let ${\mathfrak{a}}:=\fksl^{(\deftau)}(M|M)$ with 
$(\deftau,M)\in\{(1,1),(1,2),(2,2)\}$.
Let ${\mathfrak{b}}$ be the sub-Lie superalgebra of  ${\mathfrak{a}}$
generated by its Chevalley generators  $\barh_i$, $\bare_i$, $\barf_i$ ($i\in\bI$).
Assume $(\deftau,M)\in\{(1,1),(2,2)\}$. Then $\dim{\mathfrak{b}}<\infty$,
i.e., ${\mathfrak{b}}$ is a proper sub-Lie superalgebra ${\mathfrak{a}}$.
Thus we can not define $U_q({\mathfrak{a}})$ in the Drinfeld-Jimbo's way.
We can define it in a way of $RLL=LLR$.
Assume $(\deftau,M)=(1,2)$. Then ${\mathfrak{b}}={\mathfrak{a}}$ and, 
however, 
we need infinitely many defining relations to define ${\mathfrak{a}}$ by the Chevalley generators,
see \cite{IK01}, \cite{Y01}.
\end{remark}

\vspace{2cm}

Hiroyuki Yamane \par 
Department of Mathematics \par 
Faculty of Science \par  
Toyama University \par  
Gofuku, Toyama 930-8555, JAPAN \par 
e-mail: hiroyuki@sci.u-toyama.ac.jp


\begin{thebibliography}{99}
\bibitem[1]{AYY15} S.~Azam, H.~Yamane and M.~Yousofzadeh,
Classification of Finite Dimensional Irreducible Representations of
Generalized Quantum Groups via Weyl Groupoids, 
Publ.~Res.~Inst.~Math.~Sci.
51 (2015), 59--130.


\bibitem[2]{IK01} K.~Iohara and Y.~Koga, Central Extensions of Lie Superalgebras, Comment. Math. Helv. 76, 2001, 110--154.

\bibitem[3]{EG09} B.~Enriquez  and N.~Geer,
Compatibility of quantization functors of Lie bialgebras with duality and doubling operations,
Selecta Mathematica 15 (2009), 1-59.

\bibitem[4]{Y94} H.~Yamane, Quantized Enveloping Algebras Associated
with Simple Lie Superalgebras
and Their Universal R-matrices, Publ.~Res.~Inst.~Math.~Sci. 30 (1994), 15-87

\bibitem[5]{Y99} H.~Yamane, On Defining Relations of Affine Lie Superalgebras and
Affine Quantized Universal Enveloping Superalgebras, Publ.~Res.~Inst.~Math.~Sci. 35 (1999), 321-390

\bibitem[6]{Y01} H.~Yamane,
Errata to ``On Defining Relations of Affine Lie Superalgebras and Affine Quantized Universal Enveloping Superalgebras",
Publ.~Res.~Inst.~Math.~Sci. 37 (2001), 615-619


\end{thebibliography}
\end{document}